\documentclass[a4paper,11 pt]{article}

\usepackage{latexsym}
\usepackage[english]{babel}
\usepackage[utf8]{inputenc}  
\usepackage{amsmath} 
\usepackage{amssymb} 
\usepackage{mathrsfs}
\usepackage{graphicx}
\usepackage{ esint }
\usepackage[usenames]{color}
\usepackage{enumerate}
\usepackage{color}
\usepackage{dsfont}
\usepackage{bbold}
\usepackage[textsize=small]{todonotes}       

\def\epsilon{\varepsilon}

\newtheorem{lemm}{Lemma}[section]

\newtheorem{theorem}[lemm]{Theorem}

\newtheorem{ex}[lemm]{Example}
\newtheorem{defi}[lemm]{Definition}
\newtheorem{remark}[lemm]{Remark}

\def\dim{{\bf Proof. }}

\def\a1{$a_{1}$}
\def\r1{$r_{1}$}

\def\1{\{}
\def\2{\}}

\newcommand{\cvd}{\begin{flushright}$\Box$\end{flushright}}

\newcommand{\eq}{\begin{equation}}
\newcommand{\feq}{\end{equation}}
\newcommand{\be}{\begin{equation}}
\newcommand{\ee}{\end{equation}}

\newmuskip\pFqmuskip
\newcommand*\pFq[6][8]{%
  \begingroup 
  \pFqmuskip=#1mu\relax
  \mathcode`\,=\string"8000
  \begingroup\lccode`\~=`\,
  \lowercase{\endgroup\let~}\pFqcomma
  {}_{#2}F_{#3}{\left(\genfrac..{0pt}{}{#4}{#5};#6\right)}%
  \endgroup
}
\newcommand{\pFqcomma}{\mskip\pFqmuskip}

\begin{document}







\title{Self-duality of Markov processes and intertwining functions}
\author{
Chiara Franceschini
\footnote{University of Ferrara, via Macchiavelli 30, 44121 Ferrara, Italy. {\tt frnchr@unife.it}}
\and
Cristian Giardin\`a
\footnote{University of Modena and Reggio Emilia, via G. Campi 213/b, 41125 Modena, Italy. {\tt cristian.giardina@unimore.it}}
\and
Wolter Groenevelt
\footnote{Technische Universiteit Delft, DIAM, P.O. Box 5031, 2600 GA Delft, The Netherlands. {\tt w.g.m.groenevelt@tudelft.nl}}
}

\maketitle

\vspace{1.cm}

\begin{abstract}
\noindent
We present a theorem which elucidates the connection between self-duality of Markov processes and representation theory of Lie algebras. In particular, we identify sufficient conditions such that the intertwining function between two representations of a certain Lie algebra is the self-duality function of a (Markov) operator. 
In concrete terms, the two representations are associated to two operators in interwining relation.
The self-dual operator, which arise from an appropriate symmetric linear combination of them, is the generator of a Markov process.
The theorem is applied to a series of examples, including Markov processes with a discrete
state space (e.g. interacting particle systems) and Markov processes with continuous state space (e.g. diffusion processes). In the examples we use explicit representations of Lie algebras that are unitary equivalent. As a consequence, in the discrete setting self-duality functions are given by orthogonal polynomials whereas in the continuous context they are Bessel functions.
\end{abstract}

\newpage
\section{Introduction}

In the theory of interacting particle systems  \cite{Liggett, DP06}, and more generally in the theory of Markov processes,
stochastic duality plays a key role. Duality is a fundamental tool by which the analysis of the process is substantially simplified.
A  list of examples of systems that have been analyzed using  duality includes: boundary driven models of transport  and derivation of Fourier's law \cite{KMP, Spohn, GKR, CGGR2},
diffusive particle systems and their hydrodynamic limit \cite{DP06}, asymmetric interacting particle systems scaling to KPZ equation \cite{Sch97, STS, BorCorSa14, CST16},
six vertex models \cite{BorCorGor16, CorPet16}, multispecies particle models \cite{Kuan15, Kuan16, Kuan17, BelSch15}, correlation inequalities \cite{GRV10} and mathematical population genetics
\cite{M, CGGR}.

In this paper we will focus on self-duality only, which can always be thought of as a special case of duality where the dual process is an independent copy of the first one.
Indeed, self-duality provides a link between a process and its copy where two different variables play a role. The simplification of self-duality typically arises
from the fact that in the copy process only finitely many particles or variables are considered.

\medskip

There is clearly need for a deeper understanding of the origin of duality property, in particular
the sufficient conditions that would guarantee the existence of a dual process.
An algebraic approach has been proposed in a series of works \cite{GKRV, CGGR}
(one should also mention the works \cite{Schutz-Sandow94, Sch97} where a connection between stochastic duality and
symmetries of quantum spin chains was noticed).
The algebraic perspective starts from the hypothesis that the Markov generator is an element of the universal enveloping algebra of a Lie algebra. Then
the derivation of a duality relation is based on two structural ideas:
\begin{itemize}
\item[(i)] duality can be seen as a {\em change of representation}: more precisely one moves between two equivalent representations and
the intertwiner of those representations yields the duality function.

\item[(ii)] self-duality is related to reversibility of the process and the existence of symmetries, i.e.~elements that {\em commute} with the process generator.
\end{itemize}
Remarkably, this scheme can also be extended to quantum deformed algebras \cite{CGRT,CGRS16}.
The goal of this paper is to show how self-duality can be framed under the change of representation of item (i).
Recently, an independent approach has established a connection between stochastic duality and the theory of special functions.
In particular the works \cite{CED, FG17, RS17, BC17} prove that for a large class of processes duality functions are provided by
orthogonal polynomials.
This result, which has been proved following an analytic approach  -- either using structural
properties of hypergeometric polynomials \cite{FG17} or generating function methods \cite{RS17} -- has been
put in an algebraic perspective in \cite{G17}, where it is shown that orthogonal duality relations correspond to unitarily equivalent representations.
\medskip

In several examples of self-duality relations the self-duality function turns out to be symmetric in the two arguments; e.g.~in case of a polynomial duality function $p_n(x)$ it is symmetric in the degree $n$ and the variable $x$.
We wonder how the symmetry of the self-duality function could be related to stochastic self-duality.
A symmetry in the self-duality function is not an hypothesis strong enough to guarantee a self-duality relation, however there are good chances that this is the case. In addition, one would also like to see some symmetries in the generator
of the process: these two facts guarantees that a self-duality relation of the process via the symmetric function can be found.
This is the main message of this paper that will be formulated in Theorem \ref{symmetric-gen}.
As an application of the theorem we will show that several known self-duality functions can be
derived in this way and we will also derive a new self-duality relation for the so-called
Brownian momentum process \cite{GKRV}.

\subsection{Paper organization}
The rest of this paper is organized as follows, in Section \ref{two} we start by recalling the definition of self-duality and then the
main theorem of this paper is presented. Namely, we show how to relate symmetric functions with the stochastic self-duality.
{Moreover, we present under which hypothesis a function that intertwines between two Lie algebra representations is a self-duality function for an operator which is written in terms of the Lie algebra generators.
In Section \ref{three} five stochastic Markov processes are described via their generators; we will deal with three interacting particle systems and two diffusive systems.
Finally, as an application, we implement our main theorem in Section \ref{five} to prove a self-duality relation for our processes.

\section{Main results}\label{two}

All the theorems of this paper involve a self-duality relation, so we start by recalling the definition of self-duality. Since the two processes are independent copies of each other one may wonder how actually self-duality makes things simpler: thanks to a self-dual process it is possible to compute the $n-$points correlation function of the initial process
via the study of the evolution of only $n$ dual particles.
\begin{defi}[Self-duality of semigroups.]
\label{dop}
Let $X=(X_{t})_{t\geq 0}$ be a continuous time Markov process with state space $\mathcal{S}$. We say that $X$ is {\em self-dual} with duality function $D:\mathcal{S} \times \mathcal{S} \longmapsto \mathbb{R} $ if
\begin{equation}\label{du}
\mathbb{E}_{x}[D(X_{t},y)]=\mathbb{E}_{y}[D(x,Y_{t})] \ ,
\end{equation}
for all $x, y \in \mathcal{S}$ and $t\geq 0$. Here $Y$ is an independent copy of the process $X$. In \eqref{du} $ \mathbb{E}_{x} $ is the expectation w.r.t. the law of the $X$ process initialized at $ x $.
\end{defi}

Under suitable hypothesis (see \cite{JK}), the above definition is equivalent to the definition of self-duality of Markov generators.
\begin{defi}[Self-duality of generators.]
\label{dual2}
Let $L$ be a generator of the Markov process $X=(X_{t})_{t\geq 0}$.
We say that $L$ is {\em self-dual} with self-duality function $D: \mathcal{S} \times \mathcal{S} \longrightarrow \mathbb{R}$ if
\begin{equation}\label{duall}
[LD(\cdot,y)](x)=[LD(x,\cdot)](y)
\end{equation}
where we assume that both sides are well defined, i.e. $D(\cdot,y)$ and $D(x,\cdot)$ belongs to the domain of the generator.
\end{defi}
In \eqref{duall} it is understood that $L$ on the lhs acts on $D$ as a function of the first variable $x$, while $L$ on the rhs acts on $D$ as a function of the second variable $y$. In case the process is a countable Markov chain the above definition can be written in matrix notation as $LD=DL^{T}$ where $D$ is an $\mid \mathcal{S} \mid\times \mid \mathcal{S} \mid$ matrix with entries $D_{x,y} =D(x,y)$ and $L^{T}$ denotes the matrix transposition of $L$.

Definition \ref{dual2} is easier to work with, so we will always work under the assumption that the notion of self-duality is the one in equation \eqref{duall}. It is our aim to show self-duality as a change of representation, so it will be convenient to extend the definition of self-duality for operators as well.
\begin{defi} [Self-duality of operators.] \label{dual3}
Let $A$ be a generic operator with domain $\mathcal{D}(A)$. We say that $A$ is self-dual if
\begin{equation}\label{opdu}
[AD(\cdot,y)](x)=[AD(x,\cdot)](y),
\end{equation}
where $D=D(x,y)$ is the self-duality function and $D(\cdot,y), \ D(x,\cdot) \in \mathcal{D}(A)$.
\end{defi}
We are now ready to introduce the main results of this paper. Except for the Heisenberg algebra, the Lie algebras we will work with features an element in the universal enveloping algebra, the Casimir $\Omega$, which commutes with every other element of the algebra.
It is interesting to notice that, whenever the Casimir is available within the algebra, then the generator of the processes defined in Section \ref{three} can be related to it via the coproduct $\Delta$. Recall that the coproduct is defined for the Lie algebra generator $X$ as
\begin{equation}
\Delta(X)=1\otimes X + X \otimes 1
\end{equation}
and that it can be extended as an algebra morphism to the universal enveloping algebra. We will see that the generator of the process is, up to a constant, equal to the coproduct of the Casimir.

We now move on to the following theorem which is the main result of this paper. It will be applied in Section \ref{five} to explicit examples. Let $\mathcal{S}$ be a metric space, we denote by $\mathscr{F}(\mathcal{S})$ the space of real-valued functions on $\mathcal{S}$. We will also work with functions $f:\mathcal{S} \times \mathcal{S} \rightarrow \mathbb{R}$ and a linear operator $A: \mathcal{D}(A) \subset \mathscr{F}(\mathcal{S}) \rightarrow \mathscr{F}(\mathcal{S}) $.

We will need the notion of intertwining function between two operators. After its definition we present a basic example to clarify.
\begin{defi}[Intertwining function] \label{intdef}
The real-valued function $f:(x,y) \rightarrow f(x,y)$ defined on $\mathcal{S} \times \mathcal{S}$ is an intertwining function between operators $A$ and $B$ if the action of $A$ on the first variable of $f$ is equal to the action of $B$ on the second variable, i.e. $ \left( Af(\cdot,y)\right) (x)=\left( Bf(x,\cdot)\right) (y) $.
\end{defi}

If $A=B$ in definition \ref{intdef}, then $A$ is a self-dual operator (in the sense of definition \ref{dual3}) with self-duality function given by the intertwining function.

\begin{ex}
Consider the two operators acting on $g: \mathbb{R}\rightarrow \mathbb{R}$ defined as follows
\begin{equation*}
\left( A g \right) (x) := xg(x) \qquad \left( B g \right) (y) := \dfrac{\partial}{\partial y} g(y) \;.
\end{equation*}
Then $f(x,y)=e^{xy}$ in an intertwining function between $A$ and $B$ since
\[
\left( Af(\cdot,y)\right) (x) = xe^{xy} =\left( Bf(x,\cdot)\right) (y) \;.
 \]
\end{ex}

\begin{remark}[Notation] \label{notation}
Sometimes it will be convenient to have a shorter notation: if $T:\mathscr{F}(\mathcal S) \to \mathscr{F}(\mathcal S)$ is an operator and $f: \mathcal S \times \mathcal S \to \mathbb R$ a function, we write  $ T_{x}f$ for the function
\[ (x,y) \mapsto [ Tf(\cdot,y)] (x)  \]  and similarly for $ T_{y}f$. In this notation $f(x,y)$ is an intertwining function between $A$ and $B$ if $A_x f = B_y f$.
\end{remark}

For operators $A$ and $B$ we call a composition of the form
\[
S(A,B)=A^{n_{1}}B^{n_{2}}A^{n_{3}} \ldots A^{n_{k-1}}B^{n_{k}} \; , \qquad \text{for} \quad n_{1},\ldots , n_{k} \in \mathbb{N}_{0}
\]
a string in $A$ and $B$.
If $S(A,B)$ is a string of this form, then the reverse string is
\[
S^{rev}(A,B)=A^{n_{k}}B^{n_{k-1}} \ldots A^{n_{2}}B^{n_{1}}
\]
and this operation can be extended to linear combinations of strings: if
\be
C= \sum_{i=1}^{k} c_{i}S_{i}(A,B) \;,
\ee
then
\be
C^{rev}= \sum_{i=1}^{k} c_{i}S^{rev}_{i}(A,B) \;.
\ee
We are particularly interested in operators such that
$C=C^{rev}$.

\begin{theorem}[Intertwining functions, symmetries and self-duality] \leavevmode
\label{symmetric-gen}
Let $A$ and $B$ be finite order difference or differential operators on $\mathscr{F}(\mathcal{S})$, and let $f=f(x,y)$, $f:\mathcal S \times \mathcal S \to \mathbb R$ be an intertwining function between $A$ and $B$.
\begin{enumerate}
\item If $f$ is symmetric, i.e., $f(x,y)=f(y,x)$, then $f$ is an intertwining function between $B$ and $A$.

\item Suppose that $f$ is also an intertwining function between $B$ and $A$, and $C$ is a linear combination of strings in $A$ and $B$ such that $C(A,B)=C^{rev}(A,B)$.
Then $C$ is a self-dual operator with self-duality function $f$.
\end{enumerate}
\end{theorem}
\dim
For the first item, using the intertwiner hypothesis $ \left( Af(\cdot,y)\right) (x)=\left( Bf(x,\cdot)\right) (y) $ and the symmetry of $f$, namely $f(x,y)=f(y,x)$, we will show that $ \left( Af(x,\cdot)\right) (y)=\left( Bf(\cdot,y)\right) (x) $. First we show that
\be \label{uno}
[Af(\cdot,y)] (x) = [Af(y, \cdot)] (x) \;.
\ee
In the discrete case, denoting by $a_{x,x'} $ the elements of
the matrix associated to the operator $A$, we have
$$
\left( Af(y,\cdot)\right) (x) =
\sum_{x'} a_{x,x'}f(y,x')=\sum_{x'} a_{x,x'}f(x',y)=(Af(\cdot,y))(x)  \;,
$$
where we used the symmetry of the function $f$.
In the continuous case for a first order differential operator $\partial_x=\dfrac{\partial}{\partial x}$,
\[
[\partial_xf(\cdot,y)] (x) =  \lim_{h\rightarrow 0} \dfrac{f(x+h,y)+f(x,y)}{h}=\lim_{h\rightarrow 0} \dfrac{f(y,x+h)+f(y,x)}{h}=[\partial_xf(y, \cdot)] (x) \;.
\]
For a finite order differential operator $A=\sum a_{k}(x) \partial_{x_1}^{k_1} \cdots \partial_{x_L}^{k_L}$, $x=(x_1,\ldots,x_L)$, this leads to \eqref{uno}
as in the previous case, assuming $f$ is sufficiently smooth.
Our initial hypothesis that the function $f$ intertwines between the operator $A$ and
the operator $B$ implies that
\be\label{due}
\left( Af(\cdot,x)\right) (y)= \left( Bf(y,\cdot)\right) (x)
\ee
Identity \eqref{uno} holds for the operator $B$ as well, i.e.
\be\label{tre}
\left( Bf(y,\cdot)\right) (x) = \left( Bf(\cdot,y)\right) (x)\;.
\ee
Combining together (\ref{uno}), (\ref{due}) and (\ref{tre}) one proves that
 $$
 \left( Af(x,\cdot)\right) (y)=\left( Bf(\cdot,y)\right) (x) \;,
 $$
i.e., $f$ is an intertwining function for $B$ and $A$: in the notation of Remark \ref{notation} this is $B_{x}f=A_{y}f$.

For the second item observe that
 \[ (A^{n_{1}}B^{n_{2}})_{x} f = (A^{n_{1}})_{x}(A^{n_{2}})_{y} f = (A^{n_{2}})_{y}(A^{n_{1}})_{x} f=
  (A^{n_{2}}B^{n_{1}})_{y} f \;,
  \]
 assuming again that $f$ is sufficiently smooth in case of differential operators.
Iterating this procedure we get that
\be \label{string}
S(A,B)_{x}f= S^{rev}(A,B)_{y} f  \;.
\ee
So now
\[
C_{x}f = \sum_{j=1}^{k} c_{j}S_{j}(A,B)_{x}f= \sum_{j=1}^{k} c_{j}S^{rev}_{j}(A,B)_{y}f =C^{rev}_{y}f
= C_{y}f \;,  \]
where the second identity comes from \eqref{string} and the fourth identity holds due to conditions on $C$.
\cvd

\begin{ex}
Suppose, as in the previous theorem, that $f$ is an intertwining function between operators $A$ and $B$ as well as between $B$ and $A$, then examples of self-dual operators are
\begin{itemize}
\item $C_{1}=AB$.
\item $C_{2}=[A,B]^{2}=ABAB+BABA - AB^{2}A - BA^{2}B$.
\end{itemize}
\end{ex}

\begin{remark}
In our applications the operator $C$ will always be the generator of a Markov process. $C$ is the concatenation of `` building blocks '' operators that turn out to be the generators of certain Lie algebras.

\end{remark}

Theorem \ref{symmetric-gen} heavily relies on `` building blocks '' operators $A$ and $B$ so one may wonder how to construct them.
In the majority of the cases the two operators $A$ and $B$ arises naturally from the structure of the Casimir element of the underlying algebra. The next lemma shows that, whenever the generator is (in terms of) the coproduct of the Casimir, $A$ and $B$ can be found as the coproduct of two other operators.

\begin{lemm}\label{lemmasym}
If the Casimir element $\Omega$ is a linear combination of strings in $X$ and $Y$, i.e. $\Omega=h(X,Y)$, then
$\Delta(\Omega)=h(\Delta(X),\Delta(Y))$. In particular, if $\Omega=\Omega^{rev}$, then $\Delta(\Omega)=\Delta(\Omega)^{rev}$.
\end{lemm}
\dim
From the fact that the coproduct is an algebra homomorphism, it follows that the
coproduct of $\Omega$ satisfies
\[
\Delta\left( \Omega \right)=\Delta\left( h(X,Y) \right)= h(\Delta (X), \Delta (Y)).
\]
\cvd
In the applications of the next section, $h$ turns out to be a polynomial of fourth degree.
Moreover, anytime the process generator is, up to a constant, equal to the coproduct of the Casimir $L \sim \Delta(\Omega)$, it will be sufficient to look for operators $X$ and $Y$ for which the Casimir $\Omega$ is equal to $\Omega^{rev}$ instead of operators $A$ and $B$ for which the generator $L$ is equal to $L^{rev}$.
We end this section showing that, once an intertwining function between two operators is available, it can be used to find an intertwining function for the coproduct of the two operators in the following way.
\begin{lemm}\label{intertwininglemma}
If $X$ and $Y$ are two Lie algebra elements acting on $\mathscr{F}(\mathcal{S})$ and $f(x,y)$ is an intertwining function between $X$ and $Y$, then $f(x_{1},y_{1})f(x_{2},y_{2})$ intertwines $\Delta(X)$ with $\Delta(Y)$.
\end{lemm}
\dim
Using the properties of the coproduct one has
\[
[\Delta(X)f(\cdot,y_{1})f(\cdot,y_{2})](x_{1}, x_{2})=f(x_{1},y_{1})  [X f(\cdot,y_{2})] (x_{2})+ [Xf(\cdot,y_{1})](x_{1}) f(x_{2},y_{2})
\]
which, using the intertwining hypothesis, becomes
\[
f(x_{1},y_{1})  [Y f(x_{2},\cdot)] (y_{2})+ [Yf(x_{1}, \cdot)](y_{1}) f(x_{2},y_{2})=[\Delta(Y)f(x_{1},\cdot)f(x_{2},\cdot)](y_{1}, y_{2}) \;.
\]
\cvd

\section{Description of the processes}\label{three}
In this section a description of the five processes considered is given: in three of them the quantity of interest is discrete, i.e.~the number of particles, while for the other two processes presented the energy or the momentum is the continuous quantity studied.
\subsection{The Symmetric Exclusion Process}
The  Symmetric Exclusion Processes (SEP) is a family of interacting particles processes on a generic graph, labeled by the parameter $j\in\mathbb{N}/{2}$.
On the undirected and connected graph $G=(V,E)$ with $N$ sites (vertices) and edge set $E$, each site can have at most $2j$ particles, and jumps only occur when an edge exists between two sites:
jumps occur at rate proportional to the number of particles in the departure site times the number of holes in the arrival site. Note that the special case $j=1/2$ corresponds to the standard exclusion process with hard core exclusion, i.e.~each site can be either full or empty \cite{Liggett}. \\
A particle configuration is denoted by $\textbf{x}=(x_{i})_{i \in V}$ where $x_{i} \in {\{0,\ldots,2j\}}$ is interpreted as the number of particles at sites $i$. The process generator reads
\begin{equation}\label{sepj}
L^{SEP(j)}f({\bf x})=\sum_{\substack{1\le i < l \le N   \\  (i,l)\in E}} x_{i}(2j-x_{l})\left[ f({\bf x}^{i,l})-f({\bf x})\right]+(2j-x_{i})x_{l}\left[ f({\bf x}^{l,i})-f({\bf x})\right]
\end{equation}
where ${\bf x}^{i,l}$ denotes the particle configuration obtained from the configuration ${\bf x}$ by moving one particle from site $i$ to site $l$: ${\bf x}^{i,l}={\bf x}-\delta_{i}+\delta_{l}$
and $f: {\{0,1,\ldots,2j\}}^N \to \mathbb{R}$ is a function in the domain of the generator.

It is easy to verify that $L^{SEP(j)}$ conserves the total number of particles and its reversible (and thus stationary) measure is given by the homogeneous product measure with marginals the Binomial distribution
with parameters $2j>0$ and  $p\in(0,1)$, i.e.~with probability mass function
\[
\rho^{SEP}(x) = {2j \choose x} p^x (1-p)^{2j-x}\;, \qquad\qquad x\in\{0,1,\ldots,2j\}\;.
\]

\subsection{The Symmetric Inclusion Process}
The Symmetric Inclusion Processes (SIP) is a family of Markov jump processes labeled by parameter $k >0$, which can be defined in the same setting of before. In this case the state space is unbounded so that each site can have an arbitrary number of particles.
Jumps occur at rate proportional to the number of particles in the departure and the arrival sites, as the generator describes:
\begin{equation}\label{sipg}
L^{SIP(k)}f(\textbf{x})=\sum_{\substack{1\le i < l \le N   \\  (i,l)\in E}} x_{i}(2k+x_{l})\left[ f(\textbf{x}^{i,l})-f(\textbf{x})\right]+x_{l}(2k+x_{i})\left[ f(\textbf{x}^{l,i})-f(\textbf{x})\right].
\end{equation}
Detailed balance is satisfied by a product measure with marginals given by identical Negative Binomial distributions with parameters $2k>0$ and $0 < p < 1$, i.e. with probability mass function
\[
\rho^{SIP}(x) = {2k+x-1 \choose x} p^x (1-p)^{2k}\;,\qquad\qquad x\in\{0,1,\ldots\}\;.
\]

\subsection{The Brownian Energy Process} \label{BEPprocess}
The Brownian Energy Process (BEP($k$)) is a family of interacting diffusions labelled by the parameter $k$ and first introduced in \cite{GKRV} as the dual of the Symmetric Inclusion Process.\\
The generator, defined on the usual undirected connected graph $G=(V,E)$, describes kinetic energy exchange between connected sites and it reads
\begin{equation}\label{BEPgene}
L^{BEP(k)}f(\textbf{z})= \sum_{\substack{1\le i < l \le N   \\  (i,l)\in E}} \left[ z_{i}z_{j}\left( \dfrac{\partial}{\partial z_{i}} f(\textbf{z}) - \dfrac{\partial}{\partial z_{j}} f(\textbf{z})  \right)^{2}  - 2k(z_{i}-z_{j})\left( \dfrac{\partial}{\partial z_{i}} f(\textbf{z}) - \dfrac{\partial}{\partial z_{j}} f(\textbf{z})  \right)   \right].
\end{equation}
where $f: (\mathbb{R^{+}})^N \to \mathbb{R}$ is in the domain of the generator and  $\textbf{z}=(z_{i})_{i \in V}$ denotes a configuration of the process with $z_{i} \in \mathbb{R^{+}}$ interpreted as a particle energy.
It is easy to verify that the total energy of the system $ \sum_{i=1}^{N} z_{i} $ is conserved by the dynamic.\\
The stationary measure of the BEP($k$) process is given by a product of independent Gamma distribution with shape parameter $2k$ and scale parameter $\theta$, i.e.~with density function (w.r.t. Lebesgue measure)
\[
\rho^{BEP(k)}(z) = \dfrac{z^{2k-1}e^{-\frac{z}{\theta}}} {\Gamma(2k) \theta^{k}}\;.
\]

\subsection{The Brownian Momentum Process} \label{BMPprocess}
The Brownian Momentum Process (BMP) is a Markov diffusion process introduced in \cite{GK}.
On the undirected connected graph $G=(V,E)$ with $N$ vertices and edge set $E$,  the generator reads
\begin{equation}\label{BMPgene}
L^{BMP}f({\bf x})=\sum_{\substack{1\le i < l \le N   \\  (i,l)\in E}}
\left(
x_{i} \frac{\partial f}{\partial x_l}({\bf x})-   x_l \frac{\partial f}{\partial x_i}({\bf x})
\right)^{2}\;,
\end{equation}
where $f: \mathbb{R}^N \to \mathbb{R}$ is a function in the domain of the generator. A configuration is denoted by $\textbf{x}=(x_{i})_{i \in V}$ where $x_{i} \in \mathbb{R}$ has to be interpreted as a particle momentum.
A peculiarity of this process regards its conservation law: if the process is started from the configuration ${\bf x}$ then $ ||{\bf x}||_2^2 = \sum_{i=1}^{N} x^{2}_{i} $ is constant during the evolution, i.e.~the total kinetic energy is conserved.

The stationary reversible measure of the BMP process is given by a family of product measures with marginals given by independent centered Gaussian random variables with variance $\sigma^2>0$, which can be set equal to $1/2$ without loss of generality, i.e.
with density function
\[
\rho^{BMP}(x) = \dfrac{e^{-x^{2}} }{\sqrt{\pi}}
\]

\begin{remark}\label{changevar}
The BMP and BEP($1/4$) processes are related via  a \emph{change of variable}. Setting
$$
z_{i}=x_{i}^{2}
$$
one has
\be
\partial_{z_{i}}=\dfrac{1}{2x_{i}}\partial_{x_{i}}
\ee
\be
 \partial_{z_{i}}^{2}=-\dfrac{1}{4x_{i}^{2}}\partial_{x_{i}}+\dfrac{1}{4x_{i}^{2}}\partial^{2}_{x_{i}} \;.
\ee
As a consequence one finds the BEP generator $L^{BEP(k)}$ in \eqref{BEPgene} with $k=1/4$ from the BMP generator $L^{BMP}$ in \eqref{BMPgene}.
\end{remark}

\subsection{The Independent Random Walk}
The Independent Random Walkers (IRW) is one of the simplest, yet non-trivial particle system studied in the literature. It consists of independent particles that perform a symmetric continuous time random walk at rate $1$ on the undirected connected graph $G=(V,E)$. The generator is given by
\begin{equation}\label{irwgen}
L^{IRW}f(\textbf{x})= \sum_{\substack{1\le i < l \le N   \\  (i,l)\in E}}  x_{i}\left[ f(\textbf{x}^{i,l})-f(\textbf{x})\right]+x_{l}\left[ f(\textbf{x}^{l,i})-f(\textbf{x})\right]\;.
\end{equation}
The reversible invariant measure is provided by a product of Poisson distributions with parameter $\lambda>0$, i.e., with probability mass function
\[
\rho^{IRW}(x) = \dfrac{e^{-\lambda} \lambda^{x}}{x!}\,, \qquad\qquad  x\in\mathbb{N}_0\;.
\]

\section{Using our main theorem to prove stochastic self-dualities of the processes}\label{five}
This section includes five subsections where detailed examples are provided: in each subsection the natural Lie algebra and an appropriate representation on an $L^2$-space is presented. In case the representations are given in terms of unbounded operators, we assume the operators act on an appropriate dense subspace of the $L^2$-space. The duality function we encounter will be defined in terms of hypergeometric functions $\mathstrut_r F_s$.  We recall here what an hypergeometric function is.
\begin{defi}[Hypergeometric function]
The hypergeometric function is defined by the series
\[
\mathstrut_r F_s \left( {\left. \genfrac{}{}{0pt}{} {a_{1},\ldots,a_{r} } { b_{1},\ldots,b_{s} }  \right\vert {x}} \right) = \sum_{k=0}^{\infty}\dfrac{(a_{1})_{k}\cdots(a_{r})_{k}}{(b_{1})_{k}\cdots(b_{s})_{k}}\dfrac{x^{k}}{k!}
\]
where $( a )_{k}$ denotes the Pochhammer symbol defined in terms of the Gamma function as
\begin{equation*}
\left( a \right)_{k}=\dfrac{\Gamma(a+k)}{\Gamma(a)}.
\end{equation*}
\end{defi}
Whenever one of the numerator parameters, say $a_{1}$, is a negative integer $-n$, the hypergeometric function  $\mathstrut_r F_s $ is a finite sum up to $n$. In this case the hypergeometric function is a polynomial in $x$ of degree $n$, but also a polynomial in $a_j$ ($j \geq 2$) of degree $n$.

References on these are \cite{KLS} for the discrete polynomials and \cite{AAR99} for the Bessel functions.
Each sub-section ends with a theorem where the statement of a self-duality relation is proven via Theorem \ref{symmetric-gen}.

\begin{remark}[Self-duality of polynomials]
In the literature of orthogonal polynomials the symmetry of the polynomial  $p_n(x) = p_x(n)$ is often called self-duality of the polynomials.
Throughout this paper we refrain from using this name in order to not create confusion with the notion of stochastic self-duality.
\end{remark}

\begin{remark}[Working on two sites] \label{2sites}
The generators in Section \ref{three} were defined in the most general setting, i.e.~on an undirected and connected graph $G$. However, noticing that our generators only acts on two (connected) variables at a time, without loss of generality it is enough to consider the process on two sites only. The results can then be lifted to a general graph $G$ via tensor products. If on two sites the self-duality function has the form $D(\textbf{x},\textbf{y})=d(x_{1},y_{1})d(x_{2},y_{2})$, then on the graph $G$ it can be extended as
\[
D(\textbf{x},\textbf{y})= \prod_{i=1}^{|V|} d(x_{i},y_{i}) \;.
\]
\end{remark}

\subsection{$\mathfrak{su}(2)$ representation and Krawtchouk polynomials}

Generators of the $\mathfrak{su}(2)$ Lie algebra are $H$, $E$, $F$ which satisfy the following commutation relations
$$[H, E]=2E \qquad [H, F]=-2F \qquad [E, F]=H.$$
The $*$-structure is defined by $H^{ *}=H$, $E^{ *}=F$ and $F^{ *}=E$.

The Casimir element is
\[
\Omega= \frac{1}{2} H^{2} +EF+FE\;,
\]
 which is central and self-adjoint, $\Omega=\Omega^{*}$.
In this setting, $H$, $E$ and $F$ are operators on  $\mathit{l}^{2}(\mathbb{N}_{2j},\mu)$ for $j \in \mathbb{N}/2$ and $\mathbb{N}_{2j}= \lbrace 0,1, \ldots, 2j \rbrace$, where the scalar product is defined by
\begin{equation*}
\langle f,g \rangle_{L^{2}(\mu)}=\sum_{n\in \mathbb{N}_{2j} } f(n)g(n) \mu(n)
\end{equation*}
and $\mu(n)=  {2j \choose n}$.


The coproduct of the Casimir is
\begin{align*}
\Delta(\Omega)&=\Delta(\dfrac{1}{2}H^{2})+\Delta(EF) + \Delta(FE)=
\dfrac{1}{2}\Delta(H)\Delta(H)+\Delta(E)\Delta(F)+\Delta(F)\Delta(E) \\ & =
\dfrac{1}{2}(1\otimes H + H \otimes 1)^{2}+(1 \otimes E + E \otimes 1)(1 \otimes F + F \otimes 1) + (1 \otimes F + F \otimes 1)(1 \otimes E + E \otimes 1) \\& =
1\otimes \Omega + \Omega \otimes 1 + H \otimes H +2F \otimes E + 2E \otimes F
\end{align*}

On general functions $f(n)$ on $\mathbb{N}_{2j}$ the actions of the three generators are given by
\begin{align*}
\left( Hf\right) (n) & := 2(n-j)f(n)  \\
 \left( Ff \right) (n) & := (2j-n)f(n+1) \\
 \left( Ef \right) (n) & := nf(n-1)
\end{align*}
where $f(-1)=f(2j+1)=0$.
The SEP generator \eqref{sepj} restricted to two sites is
\begin{align*}
L^{SEP} & = F\otimes E+E\otimes F + \dfrac{1}{2} H \otimes H +2j^{2} \\ & =  \dfrac{1}{2} \left( \Delta (\Omega) - 1\otimes \Omega - \Omega \otimes 1 \right) + 2j^{2}  \;. \nonumber
\end{align*}

Note that $\Delta(\Omega)$ acts on a different $L^2$-space than $L^{SEP(j)}$ defined in the previous section, as discussed in Remark \ref{2sites} this can be adjusted via the tensor product.

\begin{remark}
We remark that, when acting on functions $f(n)$, the action of $ L^{SEP} $ is, up to a constant, equal to the action of
$ \Delta (\Omega) $: this follows immediately from the action of $\Omega $ itself, indeed
$\left(1\otimes \Omega f \right)(n_{1},n_{2})= 2j(j+1)f(n_{1},n_{2})$.
\end{remark}

The idea now (see \cite{KJ98} and \cite{G17}) is to look for eigenfunctions of $X_{p}= E+F-a(p)H$ with an appropriate choice of $a(p)$.
Let's start from the three term recurrence relation for the symmetric \textit{Krawtchouk polynomials}, defined via the hypergeometric function as
\begin{equation*}
K_{n}(x):=\mathstrut_2 F_1 \left( {\left. \genfrac{}{}{0pt}{} {-n,-x } { -2j }  \right\vert {\frac{1}{p}}} \right)  \qquad n,x \in \mathbb{N}_{2j}\;.
\end{equation*}
The three term recurrence relation for $ K_{n}(x) $ is
\begin{equation*}
-xK_{n}(x)=p(2j-n)K_{n+1}(x)-(2jp-2np+n)K_{n}(x)+n(1-p)K_{n-1}(x) \;.
\end{equation*}
We want to read this identity as an eigenvalue equation for $X_{p}$ with Krawtchouk polynomials as eigenfunctions. We set
\[
a(p)=\dfrac{(1-2p)}{2[p(1-p)]^{1/2}} \;,
\] so that $ k(x,n)=\left( \frac{p}{1-p}\right) ^{\frac{1}{2}(n+x)}K_{n}(x) $ is a symmetric (in $n$ and $x$) eigenfunction of the operator $X_{p}$, i.e. $$\left(X_{p}k(x,\cdot)\right) (n)=\lambda(x)k(x,n)$$ where $\lambda(x)=-\dfrac{x-j}{[p(1-p)]^{1/2}}$ is the eigenvalue.
Define
\[
H_{p}=-2[p(1-p)]^{1/2}X_{p},
\]
then $k(x,n)$ is of course also an eigenfunction of $H_p$: $(H_p k(x,\cdot))(n) = 2(x-j) k(x,n)$. Comparing this with the action of $H$, we have
\be \label{sepintert}
\left( H_p k(x,\cdot)\right)(n) = (H k(\cdot,n))(x),
\ee
i.e., $k$ is an intertwining function between $H$ and $H_p$.
We have now worked everything out in order to prove the following theorem.
\begin{theorem}
\label{SEPself-duality}
The symmetric exclusion process on two sites is self-dual with duality function $k(x_{1},n_{1})k(x_{2},n_{2})$.
\end{theorem}
\dim
The statement of the theorem follows from Theorem \ref{symmetric-gen}.
First, by Lemma \ref{intertwininglemma} $k(x_{1},n_{1})k(x_{2},n_{2})$ is an intertwining function between $\Delta(H)$ and $\Delta (H_{p})$
because of equation \eqref{sepintert}. Moreover, $k(x_{1},n_{1})k(x_{2},n_{2})$ is symmetric in $(x_{1},x_{2})$ and $(n_{1},n_{2})$, so by the first item of Theorem \ref{symmetric-gen} it is also an intertwining function between $\Delta (H_{p})$ and $\Delta(H)$.
It is left to show that $L(\Delta(H), \Delta(H_{p}))= L^{rev}(\Delta(H), \Delta(H_{p}))$ where $L$ is the generator of the SEP process.
Using Lemma \ref{lemmasym}, we can just check that $\Omega=\Omega^{rev}$ with respect to $H$ and $H_{p}$. Indeed, using the following identities
\begin{equation}\label{E and F}
F+E= X_{p}+aH \qquad F-E= \dfrac{1}{2}[X_{p},H]
\end{equation}
we have
\begin{align*}
2\Omega &= H^{2}+2EF+2FE= H^{2} + (F+E)^{2}-(F-E)^{2} = H^{2} + (X_{p}+aH)^{2} -\dfrac{1}{4}\left( \left[ X_{p}, H \right] \right)^{2} \\
& =H^{2}+\left( -\dfrac{1}{2 \sqrt{p} \sqrt{1-p}}H_{p} + aH \right)^{2} -\dfrac{1}{4}\left( \left[ -\dfrac{1}{2 \sqrt{p} \sqrt{1-p}}H_{p}, H \right] \right)^{2}.
\end{align*}
From $1+a^{2}=\frac{1}{4p(1-p)}$ we obtain
\[
\Omega= \dfrac{1}{8p(1-p)}(H^{2}+H_{p}^{2}) - \dfrac{1-2p}{8p(1-p)}(HH_{p}+H_{p}H)+\dfrac{1}{32p(1-p)}[H,H_{p}]^{2},
\]
from which we can read off that $\Omega = \Omega^{rev}$.
By Theorem \ref{symmetric-gen} the SEP generator is self-dual with $k(x_{1},n_{1})k(x_{2},n_{2})$.

%
%
\cvd


\begin{remark}[Representation]
$\mathfrak{su}(2)$ is generated by $H$ and $H_{p}$ for which we have a representation on the $n$ variable as well as a representation on the $x$ variable. Using identities \eqref{E and F} operators $E$ and $F$ can also be realised as operators on the $x$ variable, producing a different unitarily equivalent representation of the $\mathfrak{su}(2)$ algebra.
\end{remark}

%

\subsection{ $\mathfrak{su}(1,1)$ representation and Meixner polynomials}
The $\mathfrak{su}(1,1)$ algebra is the Lie algebra generated by $H$, $E$, $F$ which satisfy the following commutation relations
$$[H, E]=2E \qquad [H, F]=-2F \qquad [E, F]=H.$$
The $*$-structure is defined by $H^{*}=H$, $E^{*}=-F$ and $F^{*}=-E$.
The Casimir element is $\Omega=\frac{1}{2}H^{2}+EF+FE$ which is self-adjoint, $\Omega=\Omega^{*}$.
In this setting, $H$, $E$ and $F$ are operators on  $\mathit{l}^{2}(\mathbb{N},\mu)$ where the scalar product is defined
\begin{equation*}
\langle f,g \rangle_{L^{2}(\mu)}=\sum_{n\in \mathbb{N} } f(n)g(n) \mu(n),
\end{equation*}
where $\mu(n)=  {2k+n-1 \choose n}$ and $k>0$.
The coproduct of Casimir operator is
\begin{align*}
\Delta (\Omega) & =  1\otimes \Omega + \Omega \otimes 1 + H \otimes H +2F \otimes E + 2E \otimes F.
\end{align*}

On general functions $f(n)$ the action of the three generators is given by
\begin{align*}
\left( Hf\right) (n) & := 2(k+n)f(n)  \\
 \left( Ef \right) (n) & := (2k+n)f(n+1) \\
 \left( Ff \right) (n) & := -n f(n-1)
\end{align*}
where $f(-1)=0$.
The SIP generator \eqref{sipg} on two sites is
\begin{align*}
L^{SIP} & = -F\otimes E-E\otimes F - \dfrac{1}{2} H \otimes H +2k^{2}  \\ & = - \dfrac{1}{2} \left( \Delta (\Omega) - 1\otimes \Omega - \Omega \otimes 1 \right) + 2k^{2}  \;.
\end{align*}
Consider the symmetric \textit{Meixner polynomials }
$$ M_{n}(x):=\mathstrut_2 F_1 \left( {\left. \genfrac{}{}{0pt}{} {-n,-x } { 2k }  \right\vert {1-\frac{1}{c}}} \right)  \qquad x,n \in \mathbb{N}_{0} \;.$$
The three term recurrence relation for the Meixner polynomials is
\begin{equation*}
(c-1)xM_n(x) = c(n+2k)M_{n+1}(x) +-(n+nc+2kc)M_n(x) + nM_{n-1}(x) \;.
\end{equation*}
Let us define $X_{c}:= E-F-a(c)H$  with $a(c)=\frac{(1+c)}{2\sqrt{c}}$, for which the function $m(x,n)=c^{\frac{1}{2}(x+n)}M_{n}(x;2k,c)$ is an eigenfunction, namely

\[
\left( X_{c}m(x,\cdot) \right) (n)=\dfrac{(c-1)}{\sqrt{c}}(x+k)m(x,n) \;.
\]
Calling $H_{c}=2\dfrac{\sqrt{c}}{(c-1)}X_{c}$ we have
\[
\left(H_{c}m(x,\cdot) \right) (n)=2(x+k)m(n,x)=\left( Hm(\cdot, n)\right) (x)\;,
\]
so that $m(x,n)$ is an intertwining function between $H$ and $H_{c}$.
We have now all the ingredients to prove self-duality for the SIP process.
\begin{theorem}
\label{SIPself-duality}
The symmetric inclusion process on two sites is self-dual with duality function $m(x_{1},n_{1})m(x_{2},n_{2})$.
\end{theorem}
\dim
The proof is analogous to the proof of Theorem \ref{SEPself-duality}, note that in this case the expression for the Casimir as function of $H$ and $H_{c}$ becomes
\[
\Omega= -\dfrac{(c-1)^{2}}{8c}(H^{2}+H_{c}^{2}) + \dfrac{1-c^{2}}{8c}(HH_{c}+H_{c}H)+ \dfrac{(c-1)^{2}}{32c} [H_{c},H]^{2} \;.
\]

%
%
%
%
\cvd


%
%

\subsection{ $\mathfrak{su}(1,1)$ representation and Bessel functions}
Consider now the $\mathfrak{su}(1,1)$ Lie algebra of the previous section and the BEP($k$) process defined in Section \ref{BEPprocess}, which is a continuous diffusion, so that we will look for two continuous representations of the $\mathfrak{su}(1,1)$ algebra.
These two representations were already introduced in \cite{G17}, we recall here what we need.
Generators $H$, $E$ and $F$ are now defined on $L^{2}(\mathbb{R}^{+}, \mu_{k})$, with
\[
\mu_{k}=\dfrac{z^{2k-1}e^{-z}}{\Gamma(2k)}, \qquad k>0\;,
\]
and act on functions $f(z)$ as
\begin{align} \label{contrepre}
\left( Hf\right) (z) & :=\left(  -2z\partial_{z} - (2k-z)  \right)  f(z)\nonumber \\
 \left( Ef \right) (z) & := -\dfrac{1}{2}iz f(z) \\ \nonumber
 \left( Ff \right) (z) & := \left(  -2iz\partial_{z}^{2} -2i(2k-z)\partial_{z}+\dfrac{i}{2}(4k-z) \right) f(z)
\end{align}
where $\partial_{z}:=\dfrac{\partial}{\partial z}$. Note that in this case the $\ast$-structure is $H^{\ast}=-H$, $E^{\ast}=-E$ and $F^{\ast}=-F$ and the Casimir $\Omega$ is still self-adjoint. The BEP generator \eqref{BEPgene} on two sites is
\[
L^{BEP} = -\frac12\left( \Delta (\Omega) - 1\otimes \Omega - \Omega \otimes 1 \right) + 2k^{2}.
\]

Bessel functions of the first kind are defined in terms of hypergeometric functions as
$$ J_{\nu}(z):=\dfrac{(z/2)^{\nu}}{\Gamma(\nu +1)}\mathstrut_0 F_1 \left( {\left. \genfrac{}{}{0pt}{} { - } { \nu +1 }  \right\vert {-\frac{z^{2}}{4}}} \right) \qquad \nu > -1 \;. $$
They are solutions of the second order differential equation
\begin{equation*}
-\partial_{z}^{2} J_{\nu}(z) -\dfrac{1}{z} \partial_{z} J_{\nu}(z) +\dfrac{\nu^{2}}{z^{2}}J_{\nu}(z)= J_{\nu}(z) \;.
\end{equation*}
From the differential equation above, one infers that $J_{\nu}(zw)$ is an eigenfunction for the second order operator $T$ with eigenvalue $w^{2}$, as in \cite{G17} we have
\begin{equation} \label{operatorT}
T=-\partial_{z}^{2}-\dfrac{1}{z}\partial_{z}+\dfrac{\nu^{2}}{z^{2}} \qquad TJ_{\nu}(zw)=w^{2}J_{\nu}(zw) \;.
\end{equation}
Consider now the action of operator $F$ on the $z$ variable of $J_{\nu}(zw)$ as in equation \eqref{contrepre}, using the second order differential equation for $J_{\nu}(zw)$ in \eqref{operatorT} we can find that the eigenfunctions of $F$ are given in terms of Bessel functions $J_{2k-1}(\sqrt{zw})$, which are solutions of the following second order differential equation
\[
-2z\partial_{z}^{2} J_{2k-1}(\sqrt{zw}) -2 \partial_{z} J_{2k-1}(\sqrt{zw}) +\dfrac{(2k-1)^{2}}{2z}J_{2k-1}(\sqrt{zw})= J_{2k-1}(\sqrt{zw}) \;.
\]
Consider (\cite{G17}, Lemma 4.16) the function defined as follows,
\begin{equation}\label{besselduality}
J(z,w)=e^{\frac{1}{2}(z+w)}(zw)^{-k+\frac{1}{2}}J_{2k-1}(\sqrt{zw})\;,
\end{equation}
then $J(z,w)$ is an eigenfunction for $F$ with eigenvalue $\dfrac{1}{2}iw$, i.e.
\[
\left( FJ(\cdot, w) \right) (z)= \dfrac{1}{2}iwJ(z,w)=  \left( -E J(z, \cdot) \right) (w) \;.
\]
We see that $J(z,w)$ is an intertwining function between operators $F$ and $-E$.
\begin{theorem}
\label{BEPself-duality}
The Brownian energy process on two sites is self-dual with duality function $J(z_{1},w_{1})J(z_{2},w_{2})$.
\end{theorem}
\dim
The proof is analogous  to the proof of Theorem \ref{SEPself-duality} where the intertwined operators are $F$ and $-E$, and the Casimir is
\[
\Omega=\frac{1}{2}H^{2}+EF+FE = \frac{1}{2}[-E,F]^{2}-(-E)F-F(-E)\;,
\]
which is equal to $\Omega^{rev}$.
%
%
\cvd

%
%

\subsection{ A change of variable for Bessel functions and self-duality for the Brownian momentum process}
The idea of this section is to obtain the self-duality of the BMP process as a consequence of the change of variable highlighted in Remark \ref{changevar}. The representation \eqref{contrepre} will provide a new representation for the action of the Lie algebra generators
\begin{align} \label{bmpreo}
\left( \tilde{H}f\right) (x) & :=\left(  -x\partial_{x} - \left( \dfrac{1}{2}-x^{2} \right)   \right)  f(x) \nonumber\\
 \left( \tilde{E}f \right) (x) & := -\dfrac{1}{2}ix^{2} f(x) \\ \nonumber
 \left( \tilde{F}f \right) (x) & := \left(  -\dfrac{i}{2}\partial_{x}^{2} +ix\partial_{x}+\dfrac{i}{2}(1-x^{2}) \right) f(x) \;.
\end{align}

In equation \eqref{besselduality} we set $z=x^{2}$, $w=y^{2}$ and $k=\frac{1}{4}$ so that the candidate BMP self-duality function becomes
\[
\tilde{J}(x,y)=e^{\frac{1}{2}(x^{2}+y^{2})}|xy|^{\frac{1}{2}}J_{-1/2}(xy) = e^{\frac{1}{2}(x^{2}+y^{2})}\sqrt{\dfrac{2}{\pi }}\cos(xy) \;,
\]
where the second identity follows from the fact that, for fixed parameter $ \nu = -1/2$, Bessel functions assume the simple form of
\[
J_{-1/2}(x)=\sqrt{\dfrac{2}{\pi x}}\cos(x) \;.
\]

\begin{theorem}
\label{BMPself-duality}
The Brownian momentum process on two sites is self-dual with self-duality function
$\tilde{J}(x_{1},y_{1})\tilde{J}(x_{2},y_{2})$.
\end{theorem}
\dim
Given the $\mathfrak{su}(1,1)$ algebra representation in \eqref{bmpreo}, one could argue similarly as in Theorem \ref{BEPself-duality} to prove that $\tilde{J}(x,y)$ is indeed a self-duality function for the BMP process. However, we follow another analogous path. We show here that the self-duality for the BMP process can be obtained from the self-duality of the BEP via the change of variable in Remark \ref{changevar}.

For  the invertible operator $V: L^{2}(\mathbb{R}^{+}, \mu_{1/4})\rightarrow L^{2}_{e}(\mathbb{R}, \frac{e^{-x^{2}}}{\sqrt{\pi}})   $ given by $(Vf)(x)=f(x^{2})$  and where $L_{e}^{2}$ is the $ L^{2}$-space of even functions, we have
\begin{align*}
 \tilde{H}  & = VHV^{-1} \\
\tilde{E}  & =  VEV^{-1}  \\
\tilde{F}   & =  VFV^{-1} \;.
\end{align*}
One can now easily check that
\begin{equation} \label{generatorrelation}
L^{BMP} = VL^{BEP(1/4)}V^{-1} \;.
\end{equation}
At this point we indicate with $D^{1/4}(z,w)$ and $D(x,y)$ the self-duality functions of the BEP process with $k=1/4$ and the BMP process respectively, so that the following relation holds
\begin{equation} \label{relationford}
D(x,y) = \left( V_{x}V_{y}D^{1/4}\right) (x,y)= D^{1/4}(x^{2}, y^{2}) \;,
\end{equation}
where we use the notation of Remark \ref{notation}. For the generators this gives
\begin{align*}
L^{BMP}_{x}D& =L^{BMP}_{x}V_{x}V_{y}D^{1/4} \\
& =V_{x}L^{BEP(1/4)}_{x}V_{y}D^{1/4}   \\
& =V_{x}V_{y} L^{BEP(1/4)}_{y}D^{1/4} \\
& =V_{x} L^{BMP}_{y}V_{y}D^{1/4} \\
& =L^{BMP}_{y}D\;.
\end{align*}
Here we used that operators acting on $x$ commute with operators acting on $y$, the first and last equalities are true in virtue of equation \eqref{relationford}, the second and fourth ones both come from equation \eqref{generatorrelation}, and the third one is the self-duality of the BEP($1/4$) process in Theorem \ref{BEPself-duality}.
\cvd

%
%

\subsection{Heisenberg algebra and Charlier polynomials}
The Heisenberg algebra is the Lie algebra with generators $a$, $a^{\dagger}$ and $Z$ satisfying relations
\begin{equation*}
[a, Z]=[a^{\dagger}, Z]=0 \qquad [a, a^{\dagger}]=Z\;.
\end{equation*}
The $*$-structure is defined by $a^{ *}=a^{\dagger}$, and $Z^{ *}=Z$.
The peculiarity of this algebra is that no Casimir element is available.
In this setting, $a$, $a^{\dagger}$ and $Z$ are operators on  $\mathit{l}^{2}(\mathbb{N}_{0},\mu)$ where the scalar product is defined by
\begin{equation*}
\langle f,g \rangle_{L^{2}(\mu)}=\sum_{n\in \mathbb{N_0} } f(n)g(n) \mu(n)
\end{equation*}
and $\mu(n)= \dfrac{\lambda^{n}}{n!} $ for $\lambda > 0$.
Generators $a$, $a^{\dagger}$ and $Z$ act on functions $f(n)$ on $\mathbb{N}_{0}$ as
\begin{align} \label{heisenberggen}
(af)(n)& :=nf(n-1) \nonumber \\
(a^{\dagger}f)(n) & :=\lambda f(n+1) \\
(Zf)(n) & :=\lambda f(n) \;. \nonumber
\end{align}
In this representation the independent random walk generator \eqref{irwgen} on two sites is
\begin{equation} \label{irwgen2}
L^{IRW}=(1\otimes a - a \otimes1)(a^{\dagger}\otimes 1 -  1 \otimes a^{\dagger}) = a^{\dagger}\otimes a - 1\otimes a a^{\dagger} - a a^{\dagger}\otimes 1 + a\otimes a^{\dagger} \;.
\end{equation}
We remark that since no Casimir element is available for the Heisenberg algebra, this time we will search directly for operators for which $ L $ is equal to $L^{rev}$. To this end let us define the operator $X:=Z-a^{\dagger}$. We notice that the Heisenberg algebra is generated by $a$ and $X$,
since $Z=[a,X]$ and $a^{\dagger}=-X+[a,X]$. The generator of the IRW process in \eqref{irwgen2} becomes
\begin{equation} \label{irwgen3}
L^{IRW}=
-X\otimes a + 1\otimes a X + a X\otimes 1 - a\otimes X \;.
\end{equation}
As done before, it is time to introduce our candidate self-duality functions: the \textit{Charlier polynomials} are defined by
\[
 C_{n}(x)=\mathstrut_2 F_0 \left( {\left. \genfrac{}{}{0pt}{} {-n,-x } { - }  \right\vert {-\frac{1}{\lambda}}} \right) , \qquad x,n \in \mathbb{N}_{0}\;,
\]
and they are clearly symmetric in $x$ and $n$.
They satisfy the three term recurrence relation
\begin{equation*}
-xC_{n}(x)=\lambda C_{n+1}(x) - (n+\lambda)C_{n}(x) + nC_{n-1}(x) \;,
\end{equation*}
and the following forward shift relation
\begin{equation}\label{forward shift}
xC_{n}(x-1)=\lambda C_{n}(x) - \lambda C_{n+1}(x) \;.
\end{equation}
We conclude this section with the proof of the next theorem, by giving operators $A$ and $B$ such that the hypothesis of Theorem \ref{symmetric-gen} are satisfied.
\begin{theorem}
\label{IRWself-duality}
The independent random walk process on two sites is self-dual with self-duality function $C(x_{1},n_{1})C(x_{2},n_{2})$, where $C(x,n)=C_{n}(x)$
\end{theorem}
\dim
First, let us show that item one of Theorem \ref{symmetric-gen} is satisfied. From the definition of the Charlier polynomials we have that  $C(n,x)=C(x,n)$, so that $C(x_{1},n_{1})C(x_{2},n_{2})$ is symmetric in $(x_1,x_2)$ and $(n_1,n_2)$. Define $A=a \otimes 1 - 1 \otimes a $ and  $B=X \otimes 1 - 1 \otimes X $, then $C(x_{1},n_{1})C(x_{2},n_{2})$ is an intertwining function for $A$ and $B$. Indeed, for one site
\begin{align*}
\left( XC\left( \cdot, x \right)  \right) \left( n \right) & = \left( \left( Z- a^{\dagger} \right) C\left( \cdot, x \right)  \right) \left( n \right) \\ & = \lambda \left( C(n,x) - C(n+1,x)\right) \\ & =  xC(n,x-1) \\ &   =   \left( aC\left( n, \cdot \right) \right) \left( x\right)\;,
\end{align*}
where the second equality follows immediately from \eqref{heisenberggen} and the third one follows from equation \eqref{forward shift}.
For two sites it simply becomes
\begin{align*}
\left(  B C(x_{1}, \cdot)C(x_{2},\cdot) \right) (n_{1}, n_{2}) & = (XC(x_{1}, \cdot))(n_{1})C(x_{2},n_{2}) - C(x_{1},n_{1}) (XC(x_{2},\cdot)(n_{2}) \\ &=  (aC(\cdot,n_{1})(x_{1})C(x_{2},n_{2})  - C(x_{1},n_{1}) (XC(\cdot,n_{2}))(x_{2}) \\ &=
\left(  A C(\cdot, n_{1})C(\cdot, n_{2},) \right) (x_{1}, x_{2})  \;.
\end{align*}
For the second item one can check from \eqref{irwgen3} that the generator of the IRW process is given by $L=AB$, which is equal to $L^{rev}$.
\cvd

%


\begin{thebibliography}{empty}


\bibitem{AAR99} \label{AAR99} G. E. Andrews, R. Askey, R. Roy. {\em Special Functions}. Encycl. Math. Appl. 71, Cambridge University Press (1999).

\bibitem{BelSch15} \label{BelSch15} V. Belitsky, G.M. Schütz. Self-duality for the two-component asymmetric simple exclusion process. Journal of Mathematical Physics 56.8, 083302 (2015).


\bibitem{CED} \label{CED} C. Bernardin. Superdiffusivity of asymmetric energy model in dimensions $1$ and $2$. Journal of Mathematical Physics 49.10 103301 (2008).



\bibitem{BorCorGor16} \label{BorCorGor16} A. Borodin, I. Corwin, V. Gorin. Stochastic six-vertex model. Duke Mathematical Journal 165.3, 563--624 (2016).

\bibitem{BC17} \label{BC17} A. Borodin, I. Corwin. Dynamic ASEP, duality and continuous $ q^{-1} $-Hermite polynomials. Preprint arXiv:1705.01980 (2017).



\bibitem{BorCorSa14} \label{BorCorSa14} A. Borodin, I. Corwin, T. Sasamoto. From duality to determinants for $q$-TASEP and ASEP. The Annals of Probability 42.6, 2314--2382 (2014).


\bibitem{CGGR2} \label{CGGR2} G. Carinci, C. Giardinà, C. Giberti, F. Redig. Duality for stochastic models of transport. Journal of Statistical Physics 152.4, 657--697 (2013).

\bibitem{CGRT} \label{CGRT} G. Carinci, C. Giardinà, F. Redig, T. Sasamoto.
A generalized Asymmetric Exclusion Process with $U_q(\mathfrak{sl}_2)$ stochastic duality. Probability Theory and Related Fields, 1--47 (2014).

\bibitem{CGRS16} \label{CGRS16} G. Carinci, C. Giardinà, F. Redig, T. Sasamoto. Asymmetric Stochastic Transport Models with ${\mathscr{U}}_q (\mathfrak{su}(1, 1))$ Symmetry. Journal of Statistical Physics 163.2, 239--279 (2016).

\bibitem{CGGR} \label{CGGR} G. Carinci, C. Giardinà, C. Giberti, F. Redig. Dualities in population genetics: a fresh look with new dualities. Stochastic Processes and their Applications 125.3, 941--969 (2015).








\bibitem{CorPet16} \label{CorPet16} I. Corwin, L. Petrov. Stochastic higher spin vertex models on the line. Communications in Mathematical Physics 343.2, 651--700 (2016).

\bibitem{CST16} \label{CST16} I. Corwin, H. Shen, L-C Tsai. ASEP($q,j$) converges to the KPZ equation. Preprint arXiv:1602.01908 (2016).

\bibitem{DP06}	\label{DP06}
A. De Masi, E. Presutti. {\em Mathematical methods for hydrodynamic limits}. Springer (2006).


\bibitem{GK} \label{GK} C. Giardinà, J. Kurchan. The Fourier law in a momentum-conserving chain. Journal of Statistical Mechanics: Theory and Experiment 2005.05,  P05009 (2005).


\bibitem{GKR} \label{GKR} C. Giardinà, J. Kurchan, F. Redig. Duality and exact correlations for a model of heat conduction. Journal of Mathematical Physics 48.3, 033301 (2007).



\bibitem{FG17}\label{FG17}
C. Franceschini, C. Giardinà. Stochastic Duality and Orthogonal Polynomials. Preprint arXiv:1701.09115 (2017).


\bibitem{GKRV} \label{GKRV} C. Giardinà, J. Kurchan, F. Redig, K. Vafayi. Duality and hidden symmetries in interacting particle systems. Journal of Statistical Physics 135.1, 25--55  (2009).


\bibitem{GRV10} \label{GRV10} C. Giardinà, F. Redig, K. Vafayi. Correlation inequalities for interacting particle systems with duality. Journal of Statistical Physics 141.2, 242--263 (2010).



\bibitem{G17} \label{G17} W. Groenevelt. Orthogonal stochastic duality functions from Lie algebra representations. Preprint arXiv:1709.05997 (2017).


\bibitem{JK} \label{JK} S. Jansen, N. Kurt. On the notion(s) of duality for Markov processes, Probability Surveys 11, 59--120 (2014).

\bibitem{KMP}\label{KMP} C. Kipnis, C. Marchioro, E. Presutti. Heat flow in an exactly solvable model. Journal of Statistical Physics 27.1, 65--74 (1982).

\bibitem{KLS} \label{KLS} R. Koekoek, P.A. Lesky, R.F. Swarttouw.  \emph{ Hypergeometric Orthogonal Polynomials and their $q$-Analogues}, Springer (2010).


\bibitem{KJ98}\label{KJ98}
H.T. Koelink, J. Van der Jeugt. Convolution for orthogonal polynomials from Lie and quantum algebra representations. SIAM J. Math. Annal. 29, 794--822 (1998).

\bibitem{Kuan15}\label{Kuan15}
J. Kuan. Stochastic duality of ASEP with two particle types via symmetry of quantum groups of rank two. Journal of Physics A: Mathematical and Theoretical 49.11, 115002 (2016).

\bibitem{Kuan16}\label{Kuan16}
J. Kuan. A Multi-species ASEP($q,j$) and $q$-TAZRP with Stochastic Duality. Preprint arXiv:1605.00691 (2016).

\bibitem{Kuan17}\label{Kuan17}
J. Kuan. An algebraic construction of duality functions for the stochastic $ U_q (A_n^{(1)})$ vertex model and its degenerations. Preprint arXiv:1701.04468 (2017).


\bibitem{RS17}\label{RS17}
R. Redig, F. Sau. Duality functions and stationary product measures. Preprint arXiv:1702.07237 (2017).

\bibitem{Liggett} \label{Liggett} T. M. Liggett. \emph{Interacting particles systems}. Springer (1985).


\bibitem{M}\label{M}M. M\"ohle. The concept of duality and applications to Markov processes arising in neutral population genetics models. Bernoulli 5.5, 761--777 (1999).




\bibitem{STS} \label{STS} T. Sasamoto, H. Spohn. One-dimensional Kardar-Parisi-Zhang equation: an exact solution and its universality. Physical review letters 104.23, 230602 (2010).


\bibitem{Sch97} \label{Sch97} G.M. Schütz. Duality relations for asymmetric exclusion processes. Journal of Statistical Physics 86.5, 1265--1287 (1997).

\bibitem{Schutz-Sandow94} \label{Schutz-Sandow94} G.M. Schütz, S. Sandow. Non-Abelian symmetries of stochastic processes: Derivation of correlation functions for random-vertex models and disordered-interacting-particle systems. Physical Review E 49.4, 2726 (1994).



\bibitem{Spohn} \label{Spohn} H. Spohn. Long range correlations for stochastic lattice gases in a non-equilibrium steady state. Journal of Physics A: Mathematical and General 16.18, 4275 (1983).



\end{thebibliography}
\end{document}